\documentclass[10pt,a4paper]{article}

\usepackage[hyperfootnotes=false,pdfpage layout=SinglePage]{hyperref}
\usepackage{mathptmx}
\usepackage[scaled=.92]{helvet}
\usepackage{courier}
\usepackage{amssymb}
\usepackage{amsmath}
\usepackage{amsthm}
\usepackage{latexsym}
\usepackage{amsfonts}
\usepackage{mathrsfs}
\usepackage{setspace}

\newtheorem{theorem}{Theorem}[section]
\newtheorem{corollary}[theorem]{Corollary}
\newtheorem{lemma}[theorem]{Lemma}
\newtheorem{proposition}[theorem]{Proposition}
\theoremstyle{definition}
\newtheorem{definition}[theorem]{Definition}
\newtheorem{remark}[theorem]{Remark}
\newtheorem{example}[theorem]{Example}

\numberwithin{equation}{section}
\renewcommand{\emptyset}{\varnothing}

\newcommand{\R}{\ensuremath{\mathbb R}}    
\newcommand{\C}{\ensuremath{\mathbb C}}    
\newcommand{\N}{\ensuremath{\mathbb N}}    

\newcommand{\product}{[\cdot\,,\cdot]}
\newcommand{\hproduct}{(\cdot\,,\cdot)}

\newcommand{\calH}{\mathcal H}         \newcommand{\frakH}{\mathfrak H}
         
         \newcommand{\frakJ}{\mathfrak J}
\newcommand{\calK}{\mathcal K}         
\newcommand{\calL}{\mathcal L}         \newcommand{\frakL}{\mathfrak L}

\newcommand{\calU}{\mathcal U}

\newcommand{\la}{\lambda}
\newcommand{\veps}{\varepsilon}

\newcommand{\mat}[4]
{
   \begin{pmatrix}
      #1 & #2\\
      #3 & #4
   \end{pmatrix}
}
\newcommand{\vek}[2]
{
   \begin{pmatrix}
      #1\\
      #2
   \end{pmatrix}
}

\renewcommand{\Im}{\operatorname{Im}}
\renewcommand{\Re}{\operatorname{Re}}

\newcommand{\ran}{\operatorname{ran}}
\newcommand{\dom}{\operatorname{dom}}

\newcommand{\sap}{\sigma_{{ap}}}

\renewcommand{\sp}{\sigma_{+}}

\newcommand{\spp}{\sigma_+}
\newcommand{\smm}{\sigma_-}
\newcommand{\sess}{\sigma_{\rm ess}}

\newcommand{\Lra}{\Longrightarrow}

\newcommand{\upto}{\uparrow}
\newcommand{\downto}{\downarrow}

\newcommand{\ol}{\overline}
\newcommand{\ds}{\dotplus}
\newcommand{\wt}{\widetilde}

\newcommand{\sgn}{\operatorname{sgn}}
\newcommand{\dist}{\operatorname{dist}}

\newcommand{\slim}[1]{\underset{#1}{{\rm s}{\text -}\!\lim}\,}
\DeclareMathOperator*{\essinf}{ess\,inf}

\begin{document}
\vspace*{-.3cm}
\begin{center}
\begin{spacing}{1.7}
{\LARGE\bf Bounds on the non-real spectrum of differential operators with indefinite weights}
\end{spacing}

\vspace{1cm}
{\Large Jussi Behrndt, Friedrich Philipp and Carsten Trunk}
\end{center}

\vspace{.5cm}
\begin{abstract}\noindent
Ordinary and partial differential operators with an indefinite weight function can be viewed as bounded perturbations of
non-negative operators in Krein spaces.
Under the assumption that $0$ and $\infty$ are not singular critical points of the unperturbed operator
it is shown that a bounded additive perturbation leads to an operator 
whose non-real spectrum is contained in a compact set and with
 definite type real spectrum outside this set. The main results
 are  quantitative estimates for this set,
 which are applied to Sturm-Liouville
and second order elliptic partial differential operators with indefinite weights on unbounded domains.
\end{abstract}

\vspace{.6cm}
\noindent{\it Keywords:} Krein space, non-negative operator, definitizable operator, critical point,
non-real spectrum, spectral points of positive and negative type, indefinite weight,
singular Sturm-Liouville equation, elliptic differential operator

\vspace{.1cm}
\noindent{\it MSC 2010:} 46C20, 47B50, 34B24, 47F05, 34L15

\section{Introduction}\label{s:intro}
We consider linear
operators associated with an indefinite differential expression
\begin{equation}\label{cl}
\calL=\frac{1}{w}\,\ell,
\end{equation}
where $w\not=0$
is a real-valued, locally integrable weight function which changes its sign and $\ell$
is an ordinary or partial differential expression of the form
\begin{equation*}\label{evprob}
\ell =-\frac{d}{dx} p\frac{d}{dx}+q\qquad\text{or}\qquad \ell
=-\sum_{j,k=1}^n \frac{\partial}{\partial x_j} \,a_{jk} \,
\frac{\partial}{\partial x_k}  +a
\end{equation*}
acting on a real interval or domain $\Omega \subset \mathbb R^n$, respectively.
In the first case $p^{-1}$ and $q$ are assumed to be real
and (locally) integrable. In the second case
$a\in
L^\infty(\Omega)$ is  real and  $a_{jk}$ are $C^{\infty}$-coefficients
such that $\ell$ is formally symmetric and
uniformly elliptic on $\Omega\subset\R^n$, and,
additionally,   the weight function $w$ and its inverse are
essentially bounded.
Together with appropriate boundary conditions
 (if necessary) the differential expression $\ell$
gives rise to a selfadjoint operator $T$ in an $L^2$ Hilbert space. Multiplication with $1/w$ leads to the corresponding
indefinite differential operator
$$
A=\frac{1}{w} T
$$
associated with $\calL$ in \eqref{cl}, which is selfadjoint in a weighted $L^2$ Krein space.

Most of the existing literature for differential operators with an indefinite weight focuses on regular or left-definite problems.
The spectral properties of the operators associated to $\calL$ in the case of a regular Sturm-Liouville expression $\ell$ were
investigated in detail, we refer to \cite{cl}, the monograph \cite{z} and the detailed references therein.
Also singular left-definite Sturm-Liouville
problems are well studied. Here the selfadjoint operator $T$
associated with $\ell$ is uniformly positive and, hence, the
corresponding indefinite differential operator $A$ associated with $\calL$ in \eqref{cl}
has real spectrum with a gap around zero, cf.\ e.g.,  \cite{BE83,BB99,BBW02,KWZ01,KWZ04} and \cite{z}.
In the case $T\geq 0$ it is of particular interest whether the operator $A$ is similar to a
selfadjoint operator; necessary and sufficient similarity criteria can be
found in, e.g., \cite{kk,kkm}.
The slightly more general situation where the indefinite Sturm-Liouville operator $A$
has finitely many negative
squares or is quasi-uniformly positive
is discussed in, e.g., \cite{BKT09,bt2,CN95}.
The general non-left-definite case is more difficult,
especially the situation where the essential spectrum
of the selfadjoint operator $T$
associated with $\ell$ is no longer contained in $\mathbb R^+$.
In this case subtle problems appear, as, e.g., accumulation of non-real
eigenvalues to the real axis, see \cite{B07-a,BKT08,bp,KT09}.
The spectral properties of indefinite elliptic partial differential operators have
been investigated in, e.g., \cite{F90-1,F90,FL96,P89,P00,P09} on bounded domains and in \cite{B11,CN94} on unbounded domains.

In the non-left-definite situation the indefinite differential operator $A$
typically possesses non-real spectrum. General perturbation results for
selfadjoint operators in Krein spaces from \cite{ABT08,B07,JL79} ensure that the non-real spectrum is contained in a compact set.  
To the best of our knowledge, explicit bounds on the size of this set do not exist in the literature.
For singular indefinite Sturm-Liouville operators it is conjectured in \cite[Remark 11.4.1]{z} that the
lower bound of the spectrum of $T$ is related to a bound for the non-real eigenvalues; the numerical
examples contained in \cite{BKT08,BKT09} support this conjecture.
It is one of our main objectives to confirm this conjecture and to provide explicit
bounds on the non-real eigenvalues of indefinite Sturm-Liouville operators and indefinite elliptic partial differential operators.

For this,  we develop here  an abstract Krein space perturbation approach which is designed for applications to differential
operators with indefinite weights. Besides bounds on the non-real spectrum, our general perturbation results Theorem~\ref{t:main1} and
Theorem~\ref{t:main2}
in Section~\ref{s:main} provide also quantitative estimates for intervals containing only spectrum of positive/negative 
type in terms of the norm of the perturbation and
some resolvent integral, as well as information on the critical point $\infty$ of the perturbed operator.
The basic idea of our approach is simple:
If $T$ is a semibounded selfadjoint operator in a Hilbert space with some negative lower bound $\gamma$,
then the operator $T+\gamma$ is
non-negative or uniformly positive. Hence
$$
A_0:=A+\frac{\gamma}{w}=\frac{1}{w} (T+\gamma)
$$
(or, more generally, $A_0=G^{-1}(T+\gamma)$ with some Gram operator $G$ connecting the Hilbert and Krein space inner product)
is a non-negative operator in a Krein space and therefore the spectrum of $A_0$ is real. Moreover, the difference of $A_0$ and $A$ is a
 bounded operator. In general, a bounded
perturbation of $A_0$ may lead to unbounded non-real spectrum, but under
the additional
assumption that $0$ and $\infty$ are not singular critical points of $A_0$, the
influence of the perturbation on the non-real spectrum can be controlled. The proofs of Theorem~\ref{t:main1} and \ref{t:main2}
are based on general Krein space perturbation techniques, norm estimates and local spectral theory; they are partly inspired
by methods developed in \cite{B11,lmm}.

The abstract perturbation results from Section~\ref{s:main} are applied to ordinary and
partial differential operators with indefinite weights in Section~\ref{s:app}.
First we investigate the singular indefinite Sturm-Liouville
operator $A$ with $p=1$, a real potential $q\in L^\infty(\R)$ and the particularly simple
weight function $w(x) =\sgn(x)$. It turns out in Theorem~\ref{t:SL2} that the nonreal spectrum of $A$ is contained in the set
$$\bigl\{\lambda\in\C:\dist\,\bigl(\lambda, (-d,d)\bigr)\leq 5 \Vert q\Vert_\infty,\, |\Im\la|\le 2  \Vert q\Vert_\infty \bigr\},$$
where $-d= \,\essinf_{x\in\R}q(x)$ is assumed to be negative, and estimates on real spectral points of positive and negative type are
obtained as well.
Our second example is a
second order uniformly elliptic operator defined on an unbounded domain
$\Omega\subset \mathbb R^n$ with bounded coefficients and an essentially bounded weight function $w$ having an essentially bounded inverse.
We emphasize that the estimates for the non-real spectrum
of indefinite differential operators seem to be the first ones
in the mathematical literature.

\
\\
\textbf{Notation:}
As usual, $\mathbb C^+$ ($\mathbb C^-$) denotes the open upper
half-plane (the open lower half-plane, respectively),
$\mathbb R^+:=
(0,\infty)$ and $\mathbb R^-:=(-\infty , 0)$.  The compactification
$\mathbb R \cup \{\infty\}$
of $\mathbb R$ is denoted by $\overline{\mathbb R}$,
the compactification
$\mathbb C \cup \{\infty\}$
by $\overline{\mathbb C}$.
By $L(\calH,\calK)$ we denote the set of all bounded and everywhere defined linear operators from a Hilbert space $\calH$ to a Hilbert space $\calK$. We write $L(\calH)$ for $L(\calH,\calH)$. For a closed, densely defined linear operator $T$ in $\calH$ we denote the spectrum and the resolvent set by $\sigma(T)$ and $\rho(T)$, respectively. A point $\lambda\in\C$ belongs to the {\it approximate point spectrum} $\sap(T)$ of $T$ if there exists a sequence $(f_n)$ in $\dom T$ with $\|f_n\| = 1$ for $n\in\N$ and $(T - \la)f_n\to 0$ as $n\to\infty$. The essential
spectrum of $T$ is the set $\{\lambda \in \mathbb C : A-\lambda
\mbox{ is not Fredholm}\}$.

\section{Selfadjoint operators in Krein spaces}\label{s:saop}
Let $(\calH,\product)$ be a Krein space, let $J$ be a fixed fundamental symmetry in $\calH$ and denote by $\hproduct$
the Hilbert space scalar product induced by $J$, i.e.\ $\hproduct = [J\,\cdot\,,\cdot]$. The induced norm is denoted by $\|\cdot\|$.
 For a detailed treatment of Krein spaces and operators therein we refer to the monographs \cite{ai} and \cite{b}.

For a densely defined linear operator $A$ in $\calH$ the adjoint with respect to the Krein space inner product $[\cdot,\cdot]$ is
denoted by $A^+$. We mention that $A^+ = JA^*J$, where $A^*$ denotes the adjoint of $A$ with respect to the scalar product $\hproduct$.
The operator $A$ is called {\it selfadjoint in the Krein space $(\calH,\product)$} (or  $\product$-{\it selfadjoint}) if $A = A^+$.
This is equivalent to the selfadjointness of the operator $JA$ in the Hilbert space $(\calH,\hproduct)$.
The spectrum of a selfadjoint operator $A$ in a Krein space is in general not contained in $\R$, but its real spectral points belong to the
approximate point spectrum (see, e.g., \cite[Corollary VI.6.2]{b}),
\begin{equation*}\label{e:Rsssap}
\sigma(A)\cap\R\subset\sap(A).
\end{equation*}
A selfadjoint operator $A$ in $(\calH, \product)$ is said to be {\it non-negative} if $\rho(A)\neq\emptyset$ and if
$[Af,f]\ge 0$ holds for all $f\in\dom A$. If, for some $\gamma > 0$, $[Af,f]\geq\gamma\|f\|^2$ holds for all $f\in\dom A$,
then $A$ is called  {\it uniformly positive}. A selfadjoint operator $A$ in $(\calH, \product)$ is uniformly positive if and only
if $A$ is non-negative and $0\in\rho(A)$.

It is well-known that the spectrum of a non-negative operator $A$ in a Krein space $(\calH,\product)$ is real.
Moreover, $A$ possesses a spectral function $E_A$ on $\ol\R$ which is defined for all Borel sets
$\Delta\subset\overline\R$ whose boundary does not contain the points $0$ and $\infty$. The corresponding
spectral projection $E_A(\Delta)$ is bounded and $\product$-selfadjoint; see \cite{A79,J81,KS66,L65,L82} for further details.
The point $0$ ($\infty$) is called a {\it critical point}
of $A$ if for each Borel set $\Delta\subset\overline\R$ with $0\in\Delta$ ($\infty\in\Delta$, respectively) such that $E_A(\Delta)$ is defined the inner product $\product$ is indefinite
on the subspace $E_A(\Delta)\calH$. Moreover, if the point $0$ ($\infty$) is a critical point of $A$, it is called {\it regular} if there
exists $C > 0$ such that $\|E_A([-\veps,\veps])\|\le C$ ($\|E_A([-\veps^{-1},\veps^{-1}])\|\le C$, respectively) holds for each $\veps\in (0,1)$.
Otherwise, the critical point $0$ or $\infty$ is  called {\it singular}.

In this paper we will investigate operators which are non-negative outside of some compact set $K$ with $0\in K$,
see Definition \ref{d:lnn} below\footnote{Definition \ref{d:lnn} is slightly more general than the definition used in \cite[Definition~3.1]{BJ05}. Contrary to the definition in \cite{BJ05}, a
spectral projector for an operator
 non-negative over $\ol\C\setminus K$ in the sense of Definition  \ref{d:lnn}
corresponding to the set $\ol\C\setminus K$ does, in general, not exist.}. This notion is a generalization of the above concept
of  non-negative operators  in Krein spaces and
will be used in the study of additive bounded perturbations of non-negative operators in
Section \ref{s:main} and \ref{s:app}. For a set $\Delta\subset\C$ we denote by $\Delta^*$ the set which is obtained by  reflecting $\Delta$ in the real axis, that is, $\Delta^* = \{\ol\la : \la\in\Delta\}$.

\begin{definition}\label{d:lnn}
Let $K=K^*\subset\C$ be a compact set such that $0\in K$ and $\C^+\setminus K$ is simply connected. A selfadjoint operator
$A$ in the Krein space $(\calH,[\cdot,\cdot])$ is said to be {\em non-negative over} $\ol\C\setminus K$ if for any bounded open
neighborhood $\calU$ of $K$ in $\C$ there exists a bounded $\product$-selfadjoint projection $E_\infty$ such that with respect to the decomposition
\begin{equation}\label{e:dec_lnn}
\calH = (I-E_\infty)\calH\,[\ds]\,E_\infty\calH
\end{equation}
the operator $A$ can be written as a diagonal operator matrix
\begin{equation}\label{e:deco}
A = \mat{A_0}00{A_\infty},
\end{equation}
where $A_0$ is a bounded selfadjoint operator in the Krein space $((I-E_\infty)\calH ,\product)$
whose spectrum is contained in $\ol\calU$ and $A_\infty$ is a non-negative operator in  $(E_\infty\calH ,\product)$ with
$\calU\subset\rho(A_\infty)$.
\end{definition}

Relations \eqref{e:dec_lnn} and \eqref{e:deco} encode the non-negativity outside of the set $K$, which can also be seen in the following example.
\begin{example}
Let $A$ be the direct sum of a bounded selfadjoint operator  $A_1$ in the Krein space
$(\calH_1, \product_1)$  and of a non-negative operator  $A_2$  in the Krein space $(\calH_2, \product_2)$,
$$
A :=  A_1\times  A_2.
$$
Let $K_{r_1}$ be the closed disc around zero whose radius is the spectral radius $r_1$ of $A_1$, then for a bounded
open neighborhood $\calU$ of $K_{r_1}$ in $\C$ we set
$
E_\infty:= E_2(\ol{\mathbb R}\setminus \calU)
$, where  $E_2$ is the spectral function of $A_2$. We obtain a decomposition of the form  \eqref{e:dec_lnn}
and the operator $A$ can be written as in \eqref{e:deco},
where $A_0$  and $A_\infty$ have the properties stated in
Definition \ref{d:lnn}. Therefore, the operator $A$ is
non-negative over $\ol\C\setminus K_{r_1}$.
\end{example}

Let $A$ be non-negative over $\ol\C\setminus K$. Then from the representation \eqref{e:deco} it is seen that 
for each bounded open neighborhood $\calU$ of $K$ the non-real spectrum of $A$ is contained in $\ol\calU$. Therefore,
\begin{equation}\label{nrspec}
\sigma(A)\setminus\R\subset K.
\end{equation}
By \eqref{e:deco}
the spectral properties of $A$ and $A_\infty$ in $\ol\R\setminus K$ are the same. Therefore, we say that $\infty$ is a
{\em singular (regular) critical point} of $A$ if $\infty$ is a singular (regular, respectively) critical point of $A_\infty$.
 The following proposition is a direct consequence of \cite[Theorem~3.2]{cu} and
 the representation \eqref{e:deco}.

\begin{proposition}\label{p:curgus}
Let $A$ be a selfadjoint operator in the Krein space $(\calH,\product)$ and assume that $A$ is non-negative over $\ol\C\setminus K$.
Then $\infty$ is not a singular critical point of $A$ if and only if there exists a uniformly positive operator $W$ in the Krein
space $(\calH,\product)$ such that $W\dom A\subset\dom A$.
\end{proposition}

In order to characterize  operators which are non-negative over $\ol\C\setminus K$ in Theorem \ref{t:eqthm} below we recall
the notions of the spectral points of positive and negative type of selfadjoint operators in Krein spaces. The following definition
can be found in, e.g., \cite{j03,LcMM,lmm}.

\begin{definition}\label{d:spp}
Let $A$ be a selfadjoint operator in the Krein space $(\calH,\product)$. A point $\la\in\sap(A)$ is called a {\em spectral point of
positive {\rm (}negative{\rm )} type} of $A$ if for every sequence $(f_n)$ in $\dom A$ with $\|f_n\| = 1$ and $(A - \la)f_n\to 0$
as $n\to\infty$ we have
$$
\liminf_{n\to\infty}\,[f_n,f_n] > 0\quad\Big(\limsup_{n\to\infty}\,[f_n,f_n] < 0,\;\text{respectively}\Big).
$$
The set of all spectral points of positive {\rm (}negative{\rm )} type of $A$ will be denoted by $\spp(A)$ {\rm (}$\smm(A)$,
respectively{\rm )}. A set $\Delta\subset\C$ is said to be of {\em positive {\rm (}negative{\rm )} type} with respect to $A$
if each spectral point of $A$ in $\Delta$ is of positive type {\rm (}negative type, respectively{\rm )}.
\end{definition}

The sets $\spp(A)$ and $\smm(A)$ are contained in $\R$, open in $\sigma(A)$ and the non-real spectrum of $A$ cannot accumulate to
$\spp(A)\cup\smm(A)$. Moreover, at a spectral point $\la_0$ of positive or negative type the growth of the resolvent
is of order one in the sense of the following definition, see also \cite{ajt,j03,lmm}.

\begin{definition}\label{d:growth}
Let $A$ be a selfadjoint operator in the Krein space $\calH$
and assume that $\la_0\in\ol\R$ is not an accumulation point of $\sigma(A)\setminus\R$.
We say that the growth of the resolvent of $A$ at $\la_0$ is of order $m\ge 1$ if there exist an open neighborhood $\calU$ of $\la_0$
in $\ol\C$ and $M > 0$ such that $\calU\setminus\ol\R\subset\rho(A)$ and
\begin{equation}\label{e:growth}
\|(A - \lambda)^{-1}\| \le M\,\frac{(1+|\la|)^{2m-2}}{|\Im\la|^{m}}
\end{equation}
holds for all $\la\in\calU\setminus\ol\R$.
\end{definition}

Clearly, if $\la_0\neq\infty$, Definition \ref{d:growth} can be formulated
equivalently by replacing the enumerator in \eqref{e:growth} by $1$.
Moreover, if the growth of the resolvent of $A$ at $\la_0$ is of order $m$, then it is of order $n$ for each $n > m$.
The following theorem  gives an equivalent characterization
for non-negative operators in a neighbourhood of $\infty$.

\begin{theorem}\label{t:eqthm}
Let $A$ be a selfadjoint operator in the Krein space $\calH$ and let $K=K^*\subset\C$ be a compact set such that $0\in K$ and
$\C^+\setminus K$ is simply connected. Then $A$ is non-negative over $\ol\C\setminus K$ if and only if the following conditions are satisfied:
\begin{enumerate}
\item[{\rm (i)}]   $\sigma(A)\setminus\R\subset K$ and $(\sigma(A)\setminus K)\cap\R^\pm\subset\sigma_\pm(A)$.
\item[{\rm (ii)}]  The growth of the resolvent of $A$ at $\infty$ is of finite order.
\end{enumerate}
\end{theorem}
\begin{proof}
Let $A$ be non-negative over $\ol\C\setminus K$.
It is well-known that the growth of the resolvent of a non-negative operator in a Krein space is of finite order
for all $\lambda_0 \in \ol{\mathbb R}$ and that spectral points in $\R^+$ ($\R^-$) are of
positive (resp.\ negative) type, see, e.g., \cite[Proposition II.2.1 and Proposition II.3.1]{L82} and \cite[Section 2.1]{J88}.
Therefore, the statements (i) and (ii) are consequences of
the representation \eqref{e:deco}; cf. \eqref{nrspec}.

Conversely, assume that (i) and (ii) are satisfied. Then the operator $A$ is definitizable over $\ol\C\setminus K$, see \cite{j03}.
In particular $A$ possesses a local spectral function $E$ which is defined on all Borel sets $\Delta\subset\ol\R\setminus K$ for
which neither $\infty$ nor points of $K\cap\R$  are boundary points. For such a
set $\Delta$ the following holds:
\begin{enumerate}
\item[{\rm (a)}] $E(\Delta)$ is a bounded $\product$-selfadjoint projection and commutes
with every bounded operator which commutes with the resolvent of $A$;
\item[{\rm (b)}] $\sigma(A|E(\Delta)\calH)\,\subset\,\sigma(A)\cap\ol\Delta$;
\item[{\rm (c)}] $\sigma(A|(I - E(\Delta))\calH)\,\subset\,\sigma(A)\setminus\operatorname{int}(\Delta)$, where
$\operatorname{int}(\Delta)$ is the interior of $\Delta$ with respect to the topology of
$\ol\R$;
\item[{\rm (d)}] If, in addition, $\Delta$ is a neighbourhood of $\infty$
(with respect to the topology of
$\ol\R$), then $A|(I-E(\Delta)) \calH$ is a bounded operator;
\end{enumerate}
cf. \cite[Definition 3.13, Theorems 3.15 and 4.8, Remark 4.9]{j03} and note that (d) follows from the definition of the local spectral function
\cite[Definition 3.13]{j03} and the definition of the extended spectrum $\wt\sigma(A)$
in \cite{j03}.

Let $\calU$ be a bounded open neighborhood of $K$ in $\C$ and set $E_\infty := E(\ol\R\setminus\calU)$. By (a), with respect to the decomposition
$$
\calH = (I - E_\infty)\calH\,[\ds]\,E_\infty\calH,
$$
the operator $A$ is represented by
$$
A = \mat{A_0}00{A_\infty},
$$
and by (b) the spectrum of the operator $A_\infty$ is real, $\sigma(A_\infty) \subset  \mathbb R \setminus \cal U$. According to (i) and (ii),
the operator $A_\infty$ is definitizable over $\ol\C$ and thus definitizable, see \cite[Theorem 4.7]{j03}.
Therefore, \cite[Corollary 3 of Proposition II.5.2]{L82} implies that $A_\infty^{-1}$ (and thus $A_\infty$) is non-negative in $E_\infty\calH$.
By (d), the operator $A_0= A|(I-E_\infty) \calH$ is bounded and, by (c),
$\sigma(A_0) \subset \sigma(A) \setminus \operatorname{int}(\ol\R \setminus \calU)\subset \ol\calU$.\qed
\end{proof}

\section{Bounded selfadjoint perturbations of non-negative operators}\label{s:main}
In this section we prove two abstract results on additive bounded perturbations of non-negative (and some closely connected class of)
 operators
in Krein spaces that lead to perturbed operators which are non-negative over some neighborhood of
infinity. In both cases the neigh\-borhood is given in quantitative terms. The results will be applied in Section \ref{s:app} to
singular indefinite Sturm-Liouville operators and to second order elliptic operators with indefinite weights.

\subsection{Two perturbation results}\label{ss:pert}
The following notation will be useful when formulating our main results below: For a set $\Delta\subset\R$ and $r > 0$ we define
$$
 K_r(\Delta) := \{z\in\C : \dist(z,\Delta)\,\le\,r\}.
$$
Our first main theorem concerns bounded selfadjoint perturbations of non-negative operators in Krein spaces.

\begin{theorem}\label{t:main1}
Let $(\calH,\product)$ be a Krein space with fundamental symmetry $J$ and norm
$\|\cdot\| = [J\cdot,\cdot]^{1/2}$. Let $A_0$ be a non-negative operator in $(\calH,\product)$ such that $0$ and
$\infty$ are not singular critical points of $A_0$, and $0\notin\sigma_p(A_0)$. Furthermore, let $V\in L(\calH)$
be a selfadjoint operator in $(\calH,\product)$. Then the following holds:
\begin{itemize}
 \item [{\rm (i)}] If $V$ is non-negative, then $A_0+V$ is a
  non-negative operator in $(\calH,\product)$.
 \item [{\rm (ii)}] If $V$ is not non-negative, then $A_0+V$ is non-negative over
$
\ol\C\setminus K_r((-d,d)),
$
where
\begin{equation}\label{e:estimates1}
r = \frac{1+\tau_0}{2}\|V\|,\qquad\quad d = -\frac{1+\tau_0}{2}\,\min\sigma(JV),
\end{equation}
and
\begin{equation}\label{tau0}
\tau_0 = \frac 1 \pi\,\limsup_{n\to\infty}\,\left\|\int_{1/n}^n\,\left((A_0 + it)^{-1} + (A_0 - it)^{-1}\right)\,dt\,\right\|<\infty.
\end{equation}
\end{itemize}
Moreover, in both cases  $\infty$ is not a singular critical point of $A_0 + V$.
\end{theorem}

The following simple example shows that without further assumptions a bounded selfadjoint 
perturbation $V$ of a non-negative or uniformly positive $A_0$ in $(\calH,\product)$ may lead to an operator $A_0+V$ with
unbounded non-real spectrum.

\begin{example}
Let $(\cal K,(\cdot,\cdot))$ be a Hilbert space and let $H$ be an unbounded selfadjoint operator in $\cal K$ such that $\sigma(H)\subset 
(0,\infty)$. Equip $\cal H:=\cal K\oplus\cal K$ with the Krein space inner product 
\begin{equation*}
\left[\begin{pmatrix}
       k_1 \\ k_2
      \end{pmatrix},\begin{pmatrix}
       l_1 \\ l_2
      \end{pmatrix}\right] :=(k_1,l_2)+(k_2,l_1),\qquad \begin{pmatrix}
       k_1 \\ k_2
      \end{pmatrix},\begin{pmatrix}
       l_1 \\ l_2
      \end{pmatrix}\in\cal H,
\end{equation*}
and consider the operators 
\begin{equation*}
 A_0=\begin{pmatrix} 0 & I \\ H & 0\end{pmatrix},\quad 
V=\begin{pmatrix} 0 & -I  \\ 0 & 0\end{pmatrix},\quad\text{and}\quad
A_0+V =\begin{pmatrix} 0 & 0 \\ H & 0\end{pmatrix}.
\end{equation*}
It is easy to see that $A_0$ is a non-negative operator and $V$ is a bounded selfadjoint operator 
in the Krein space $(\calH,\product)$. Moreover, as $\dom(A_0+V)=\dom H\oplus{\cal H}$ we conclude $\ran(A_0+V-\lambda)\not=\cal H$ 
for every $\lambda\in\C$, that is, $\sigma(A_0+V)=\C$.
\end{example}

Our second main result applies to operators of the form $A = G^{-1}T$, where
$T$ is a semibounded selfadjoint operator in a Hilbert space $(\calH,\hproduct)$ and $G$ is a selfadjoint bounded and boundedly invertible
operator in the Hilbert space $(\calH,\hproduct)$. In other words, $A$ is selfadjoint in the Krein space $(\calH,(G\cdot,\cdot))$
and $A + \eta G^{-1}$ is uniformly positive for suitable $\eta$. Such a situation arises, e.g., when considering
elliptic differential operators with an indefinite weight function; cf. Section 4.2.

\begin{theorem}\label{t:main2}
Let $(\calH,\hproduct)$ be a Hilbert space, $\|\cdot\| = \hproduct^{1/2}$ and $\product=(G\cdot,\cdot)$ be as above.
Let $A$ be a selfadjoint operator in the Krein space $(\calH,\product)$ such that for some $\gamma > 0$ we have
$$
[Af,f]\ge -\gamma\|f\|^2\quad\text{for all }f\in\dom A.
$$
Assume furthermore that for some $\eta > \gamma$
\begin{equation*}\label{e:tau_nu}
\tau_\eta := \frac 1 \pi\,\limsup_{n\to\infty}\,\left\|\int_{-n}^n\,(A + \eta G^{-1} - it)^{-1}\,dt\,\right\| < \infty,
\end{equation*}
holds. Then $A$ is non-negative over $\ol\C\setminus K_r((-r,r))$, where
\begin{equation*}\label{e:estimates2}
r = \eta\,\frac{1+\tau_\eta}{2}\,\|G^{-1}\|.
\end{equation*}
Moreover, $\infty$ is not a singular critical point of $A$.
\end{theorem}

\begin{remark}\label{Mai2012}
Setting $A_0 := A + \eta G^{-1}$ in Theorem \ref{t:main2} and $V:=-\eta G^{-1}$ we have $A = A_0 + V$ and
hence Theorem~\ref{t:main2} can also be seen as a variant of Theorem~\ref{t:main1}. Here
$A_0=A + \eta G^{-1}$ is uniformly positive so that $0\in\rho(A_0)$ is not a critical point of $A_0$ and
the entire imaginary axis belongs to $\rho(A + \eta G^{-1})$. It follows
from \cite[Theorem 1.2]{v} and \cite[Lemma 1]{j82} that $\tau_\eta < \infty$ is equivalent to $\infty$ not being
a singular critical point of $A_0$.
\end{remark}

In what follows we describe  the structure of the proofs of Theorem~\ref{t:main1}~(ii) and Theorem~\ref{t:main2}.
Let $E$ denote the spectral function of the non-negative operator $A_0$. As $0$ and $\infty$ are not
singular critical points of $A_0$, the spectral projections $E_+ := E((0,\infty))$ and $E_- := E((-\infty,0))$, and
the corresponding spectral subspaces $\calH_\pm := E_\pm\calH$ of $A_0$ exist. From $0\notin\sigma_p(A_0)$ it follows that
\begin{equation}\label{e:dec_expl}
\calH = \calH_+\,[\ds]\,\calH_-
\end{equation}
 and $(\calH_\pm,\pm\product)$ are Hilbert spaces, see, e.g.,
 \cite[Proposition 5.2 and Theorem 5.7]{L82}. Therefore
 \eqref{e:dec_expl} is a fundamental decomposition of the Krein space
$(\calH,\product)$. With respect to this decomposition the operator $A_0$ and the fundamental symmetry $\wt J$ corresponding to
\eqref{e:dec_expl} can be written as operator matrices:
\begin{equation}\label{e:AJ}
A_0 = \mat{A_{0,+}}00{A_{0,-}}\quad\text{and}\quad\wt J = \mat I 0 0 {-I}.
\end{equation}
Note that the operator $\pm A_{0,\pm}$ is a selfadjoint non-negative operator in the Hilbert space $(\calH_\pm,\pm\product)$. Hence,
$A_0$ is selfadjoint in the Hilbert space $(\calH,\hproduct_\sim)$, where \label{p:idea}
\begin{equation}\label{e:sim}
(f,g)_\sim := [\wt Jf,g] = [E_+f,g] - [E_-f,g],\qquad f,g\in\calH.
\end{equation}
Now, let $V\in L(\calH)$ be a bounded selfadjoint operator in the Krein space $(\calH,\product)$. Then, with respect to
the decomposition \eqref{e:dec_expl} it admits an operator matrix representation
\begin{equation}\label{e:V}
V = \mat{V_+}{B}{C}{V_-}.
\end{equation}
From the selfadjointness of $V$ in $(\calH,\product)$ one concludes that $V_\pm$ is selfadjoint in the Hilbert space $(\calH_\pm,\pm\product) = (\calH_\pm,\hproduct_\sim)$ and
\begin{equation}\label{e:sim_adjoint}
C = -B^{\wt *},
\end{equation}
where $B^{\wt *}$ denotes the adjoint of $B$ with respect to the scalar product $\hproduct_\sim$. Hence, the perturbed operator $A := A_0 + V$ is represented by
\begin{equation}\label{e:om}
A = \mat{A_{0,+} + V_+}{B}{-B^{\wt *}}{A_{0,-} + V_-}.
\end{equation}
For operators as in \eqref{e:om} we show in Theorem \ref{l:lmm}
below that
 \begin{itemize}
 \item $\R\setminus K_{\|B\|_\sim} (\sigma(A_{0,-} + V_-))$ is of positive type,
 \item $\R\setminus K_{\|B\|_\sim} (\sigma(A_{0,+} + V_+))$  is of negative type,
 \item the non-real spectrum of $A$ is contained in
\begin{equation}\label{BadHersfeld}
K:=K_{\|B\|_\sim}\bigl(\sigma(A_{0,+} + V_+)\bigr)\cap
K_{\|B\|_\sim} \bigl(\sigma(A_{0,-} + V_-)\bigr),
\end{equation}
 \item and the growth of the resolvent of $A$ at $\infty$ is of order one.
 \end{itemize}
Here $\|\cdot\|_\sim$ stands for the operator norm induced by \eqref{e:sim}.
Then, by Theorem \ref{t:eqthm}, the operator $A$ is non-negative over
$\ol\C\setminus K$, where $K$ is as in \eqref{BadHersfeld}.
In a final step it remains to bound the quantities
$$
\|B\|_\sim \quad \mbox{and} \quad  K_{\|B\|_\sim} (\sigma(A_{0,\pm} + V_\pm))
$$
 by the quantities
$r$ and $\tau_0$ as in the statement of Theorem \ref{t:main1}.

\subsection{Some auxiliary statements}\label{ss:aux}
In this section we formulate and prove some
auxiliary statements which will be used in the proofs Theorems~\ref{t:main1} and \ref{t:main2}.
For the special case of bounded operators the first part of the following theorem was already proved in
\cite[Theorem 4.2]{lmm}, see also \cite{LLMT05} and \cite[Proposition~2.6.8]{T}, \cite[Theorem 5.5]{T09} for similar results with 
unbounded entries on the diagonal.

\begin{theorem}\label{l:lmm}
Let $S_+$ and $S_-$ be selfadjoint operators in some Hilbert spaces $\frakH_+$ and $\frakH_-$, respectively,
let $M\in L(\frakH_-,\frakH_+)$ and define the operators
$$
S := \mat{S_+}{M}{-M^*}{S_-}\quad\text{and}\quad \frakJ := \mat{I}00{-I}
$$
in the Hilbert space $\frakH := \frakH_+\oplus\frakH_-$. Then the operator $S$ is selfadjoint in the Krein space
$(\frakH,(\frakJ\cdot,\cdot))$ and with $\nu := \|M\|$ the following statements hold:
\begin{enumerate}
\item[{\rm (i)}]   $\sigma(S)\setminus\R\,\subset\, K_\nu(\sigma(S_+))\cap K_\nu(\sigma(S_-))$.
\item[{\rm (ii)}]  $\R\setminus K_\nu(\sigma(S_-))$ is of positive type with respect to $S$.
\item[{\rm (iii)}] $\R\setminus K_\nu(\sigma(S_+))$ is of negative type with respect to $S$.
\end{enumerate}
Moreover, if $S_+$ is bounded from below and $S_-$ is bounded from above {\rm (}or vice versa{\rm )}, the non-real
spectrum of $S$ is bounded and the growth of the resolvent of $S$ at $\infty$ is of order one.
\end{theorem}

\begin{proof}
We set $\product := (\frakJ\cdot,\cdot)$. Let $\la\in \mathbb C \setminus K_\nu(\sigma(S_-))$ and
$\alpha := \nu/\dist(\la,\sigma(S_-)) < 1$.
We claim that for some $r(\alpha)>0$ and all $f\in\dom S$
the following implication
\begin{equation}\label{Februar2012}
\|(S - \la)f\|\le\frac{1-\alpha^2}{4\alpha}\nu\|f\|\quad\Lra\quad [f,f]\,\ge\,r(\alpha)\|f\|^2
\end{equation}
holds. In fact, set $\veps := \frac{1-\alpha^2}{4\alpha}$ and let $f\in\dom S$ with
$$
\|(S - \la)f\|\le\veps \nu\|f\|.
$$
With respect to the decomposition $\frakH = \frakH_+\oplus\frakH_-$ we write $f$ and $g := (S - \la)f$ as column vectors
$$
f = \vek{f_+}{f_-}\in\dom S_+\oplus\dom S_-, \quad g = \vek{g_+}{g_-}.
$$
Then $g_- = -M^*f_+ + (S_- - \la)f_-$, or, equivalently,
$$
f_- = (S_- - \la)^{-1}g_- + (S_- - \la)^{-1}M^*f_+.
$$
As $\|g_-\|\le\|g\| = \|(S - \la)f\|\le\veps \nu\|f\|$ and $\|(S_- - \la)^{-1}\| = \alpha/\nu$ this yields
$$
\|f_-\|\,\le\,\alpha\veps\|f\| + \alpha\|f_+\|
$$
and hence
$
\|f_-\|^2\,\le\,\alpha^2\big(\veps^2\|f\|^2 + 2\veps\|f\|^2 + \|f_+\|^2\big).
$
Using $\|f_+\|^2 = \|f\|^2 - \|f_-\|^2$ we conclude
$$
\|f_-\|^2\,\le\,\frac{\alpha^2}{1+\alpha^2}(1 + \veps)^2\|f\|^2.
$$
Hence, as $[f,f] = \|f_+\|^2 - \|f_-\|^2 = \|f\|^2 - 2\|f_-\|^2$, we obtain $[f,f]\,\ge\,r(\alpha)\|f\|^2$, where
\begin{equation}\label{qay}
r(\alpha) := 1 - \frac{2\alpha^2}{1+\alpha^2}(1 + \veps)^2 = \frac{(1 - \alpha)^2(1 + \alpha)(7 - \alpha)}{8(1+\alpha^2)} > 0,
\end{equation}
i.e., the implication \eqref{Februar2012} holds.
It follows that $(\C\setminus K_\nu(\sigma(S_-)))\cap\sap(S)\subset\sp(S)$. Hence (ii) is proved  and, as $\sp(S)$ is real, the non-real
spectrum of $S$ satisfies
\begin{equation}\label{e:minus_part}
\sigma(S)\setminus\R\,\subset\, K_\nu(\sigma(S_-)).
\end{equation}

Similarly, as in the proof of \eqref{Februar2012}, one proves that for
 $\la\in \mathbb C \setminus K_\nu(\sigma(S_+))$,  $f\in\dom S$, and for $\beta := \nu/\dist(\la,\sigma(S_+))<1$ the implication
$$
\|(S - \la)f\|\le\frac{1-\beta^2}{4\beta}\nu\|f\|\quad\Lra\quad [f,f]\,\le\,-r(\beta)\|f\|^2
$$
holds with $r(\beta)>0$ as in \eqref{qay}. This shows (iii) and $\sigma(S)\setminus\R\,\subset\, K_\nu(\sigma(S_+))$.
Together with \eqref{e:minus_part} we obtain (i).

Assume now that $S_+$ is
bounded from below and $S_-$ is bounded from above. Then the non-real spectrum of $S$ is bounded by (i) and it remains to show
that the growth of the resolvent of $S$ is of order one at $\infty$. Note first that for
$$
S_0=\begin{pmatrix}
     S_+ & 0 \\ 0 & S_-
    \end{pmatrix}
\qquad\text{and}\qquad
V=\begin{pmatrix}
   0 & M \\ -M^* & 0
  \end{pmatrix},
$$
and $\lambda\in\C$ such that $\vert\Im\lambda\vert\geq 2\Vert V\Vert = 2\nu$ we have $\Vert (S_0-\lambda)^{-1}V\Vert\leq \frac{1}{2}$ and
\begin{equation}\label{abs1}
\Vert(S-\lambda)^{-1}\Vert = \bigl\Vert\bigl((S_0-\lambda)(I+(S_0-\lambda)^{-1}V)\bigr)^{-1}\bigr\Vert\leq\frac{2}{\vert\Im\lambda\vert}.
\end{equation}
For $\lambda\in\C$ with $\vert\Im\lambda\vert < 2\nu$ and $\dist(\la,\sigma(S_-)) > \nu + \delta$ with some fixed $\delta > 0$ again define the value $\alpha = \nu/\dist(\la,\sigma(S_-))$. Then $\alpha\in (0,\delta')$, where $\delta' = \nu/(\nu + \delta)\in (0,1)$, and for $f\in\dom S$ we either have
\begin{equation}\label{abs2}
\|(S - \la)f\|\ge\frac{1-\alpha^2}{4\alpha}\nu\|f\| > \frac{1-\alpha^2}{8\alpha}\vert\Im\lambda\vert\,\|f\|,
\end{equation}
or, by \eqref{Februar2012},
$$
r(\alpha)|\Im\la|\|f\|^2\le|\Im [\la f,f]| = |\Im [(S - \la)f,f]|\le\|(S - \la)f\|\|f\|,
$$
and hence
\begin{equation}\label{abs3}
\|(S - \la)f\|\geq r(\alpha)|\Im\la|\,\|f\|.
\end{equation}
For $\la\in\C$ with $\vert\Im\lambda\vert<2\nu$ and $\dist(\la,\sigma(S_+)) > \nu + \delta$ the estimates \eqref{abs2} and \eqref{abs3} hold in a similar form. Therefore, \eqref{abs1}-\eqref{abs3} imply that for all non-real $\la\in\C$ with $\dist(\la,\sigma(S_-)) > \nu + \delta$ or $\dist(\la,\sigma(S_+)) > \nu + \delta$ we have
$$
\Vert(S-\lambda)^{-1}\Vert\leq\frac{C}{\vert\Im\lambda\vert}
$$
with some $C > 0$ which does not depend on $\la$. This shows that the growth of the resolvent $(S-\lambda)^{-1}$ is of order one at $\infty$.\qed
\end{proof}

\begin{remark}\label{r:real}
If $\dist(\sigma(S_+),\sigma(S_-))\ge 2\|M\|$, then the spectrum of $S$ in Theorem~\ref{l:lmm} is real.
This result can be improved in certain special cases, cf.\ \cite{AMS09} and also \cite{AMT10}, where sharp norm bounds on the operator
angles between reducing subspaces of $S$ and $S_+\oplus S_-$ are given.
We mention here also \cite{j93,j00} 
for comparable results in the study of operators of Klein-Gordon type,
where bounded perturbations of  operators with unbounded off-diagonal
entries are investigated.
\end{remark}

The following simple example shows that the bounds on the non-real spectrum of $S$ in Theorem \ref{l:lmm} are sharp.
\begin{example}
Let $\frakH_+ = \frakH_- = \C$ and $S_\pm = \pm 1$, $M = z\in\C$ and hence
$$
S = \mat{1}{z}{-\ol z}{-1}.
$$
Then $\dist(\sigma(S_+),\sigma(S_-)) = 2$. If $|z|\le 1$,
then $\sigma(S) = \{\pm\sqrt{1 - |z|^2}\}$ is real;
cf.\ Remark~\ref{r:real}. If $|z| > 1$, then $\sigma(S) = \{\pm i\sqrt{|z|^2 - 1}\}$ and in this case the eigenvalues of $S$ lie on the boundary
of $ K_{|z|}(\{1\})\cap K_{|z|}(\{-1\})$.
\end{example}

\begin{lemma}\label{l:veselic}
Let $T$ be a selfadjoint operator in the Hilbert space $(\frakH,\hproduct)$ with $0\notin\sigma_p(T)$ and let $E$ be its spectral function.
If $(a_n)$ and $(b_n)$ are sequences in $[0,\infty)$ such that $a_n\downto 0$ and $b_n\upto\infty$ as $n\to\infty$, and
$0\in\rho(T)$ if $a_k=0$ for some $k\in\N$, then
\begin{equation*}\label{e:strong_lim}
E(\R^+) - E(\R^-) = \frac 1 \pi\;\,\slim{n\to\infty}\!\int_{a_n}^{b_n}\,\left((T + it)^{-1} + (T - it)^{-1}\right)\,dt\,.
\end{equation*}
\end{lemma}
\begin{proof}
First of all we observe that
$$
(T + it)^{-1} + (T - it)^{-1} = 2T(T^2 + t^2)^{-1}.
$$
For all $f,g\in\frakH$ we have
\begin{align*}
\int_{a_n}^{b_n}\,\big(T(T^2 + t^2)^{-1}f,g\big)\,dt
&= \int_{a_n}^{b_n}\,\int_{\R}\,\frac s {s^2 + t^2}\,d(E_sf,g)\,dt\\
&= \int_{\R}\,\int_{a_n}^{b_n}\,\frac s {s^2 + t^2}\,dt\,d(E_sf,g)\\
&= \int_{\R}\,\big(\arctan(b_n/s) - \arctan(a_n/s)\big)\,d(E_sf,g).\\
&= \big(\arctan(b_nT^{-1})f - \arctan(a_nT^{-1})f,g\big).
\end{align*}
Fubini's theorem can be applied here since on $\R\times [a_n,b_n]$ the integrand is bounded and the measure $d(E_sf,g)\otimes dt$
considered as a measure in $\mathbb R \times [a_n,b_n]$ is finite. Hence,
$$
\int_{a_n}^{b_n}\,2T(T^2 + t^2)^{-1}f\,dt = 2\arctan(b_nT^{-1})f - 2\arctan(a_nT^{-1})f.
$$
Now, we apply \cite[Theorem VIII.5]{rs1} and observe that $2\arctan(a_nT^{-1})f$ tends to zero and $2\arctan(b_nT^{-1})f$ tends to
$$
\pi\sgn(T^{-1})f = \pi\sgn(T)f = \pi(E(\R^+) - E(\R^-))f
$$
as $n\to\infty$.\qed
\end{proof}

A subspace $\frakL$ of a Krein space $(\frakH,\product)$ is called {\it uniformly positive} if  there exists $\delta > 0$ such that
$$
[f,f]\,\ge\,\delta\|f\|^2 \quad \mbox{for all } f\in\frakL.
$$
Here $\|\cdot\|$ is a Hilbert space norm  on $\frakH$
with respect to which $\product$ is continuous.
All such norms are equivalent, see
\cite[Proposition I.1.2]{L82}.

\begin{lemma}\label{l:theoretic}
Let $G$ be a bounded and boundedly invertible selfadjoint operator in a Hilbert space $(\frakH,\hproduct)$
and set $\product := (G\cdot,\cdot)$. Let $\frakL$ be a closed uniformly positive subspace in the Krein space $(\frakH,\product)$.
If $P$ {\rm (}$E${\rm )} denotes the orthogonal projection in the Hilbert space $(\frakH,\hproduct)$ {\rm (}in the Krein space $(\frakH,\product)$,
respectively{\rm )} onto $\frakL$, then $P(G|\frakL)$ is a uniformly positive selfadjoint operator in $(\frakL,\hproduct)$,
and its inverse is given by
$$
\bigl(P(G|\frakL)\bigr)^{-1} = E\bigl(G^{-1}|\frakL\bigr).
$$
\end{lemma}
\begin{proof}
As $\frakL$ is uniformly positive, for all $f\in\frakL$ we have
$$
\big(P(G|\frakL)f,f\big) = (Gf,f) = [f,f]\,\ge\,\delta\|f\|^2
$$
with some $\delta > 0$. Thus, $P(G|\frakL)$ is a uniformly positive operator in $(\frakL,\hproduct)$ and, in particular,
$0\in\rho(P(G|\frakL))$. Let $g\in\frakL$
 and set $f := E(G^{-1}|\frakL)g\in\frakL$. Since for all $h\in\frakL$
$$
(g - Gf,h) = \bigl(G(I - E)G^{-1}g,h\bigr) = \bigl[(I - E)G^{-1}g,h\bigr] = 0,
$$
we conclude $g - Gf\in\frakL^\perp$ and hence $g = PGf = P(G|\frakL)f$. Therefore,
$$
\bigl(P(G|\frakL)\bigr)^{-1}g = f = E\bigl(G^{-1}|\frakL\bigr)g,
$$
which proves the lemma.\qed
\end{proof}

\subsection{Proof of Theorem~\ref{t:main1}}\label{ss:proof_main1}
Assertion (i) in Theorem~\ref{t:main1} is simple: As both $A_0$ and $V$ are non-negative in the Krein space $(\calH,\product)$ we have
$$[(A_0+V)f,f]\geq 0\quad\text{for all}\quad f\in\dom (A_0+V)=\dom A_0.$$
The operator $A_0$ is selfadjoint in the Hilbert space  $(\calH,\hproduct_\sim)$ (cf.\ \eqref{e:AJ} and \eqref{e:sim}) and,
as $V$ is bounded,
$\rho(A_0+V)$ is nonempty. Thus $A_0+V$ is a non-negative operator
in $(\calH,\product)$ and (i) in Theorem~\ref{t:main1} is proved.

Let $E$ be the spectral function of the non-negative operator $A_0$, define the projections $E_+ := E((0,\infty))$, $E_-:=E((-\infty,0))$ and the subspaces $\calH_\pm := E_\pm\calH$. Then with respect to the fundamental decomposition \eqref{e:dec_expl} the operators $A_0$ and $V$
have the form in \eqref{e:AJ} and \eqref{e:V}, respectively. The fundamental symmetry $\wt J$ associated with the decomposition
\eqref{e:dec_expl} is also given in \eqref{e:AJ}. The inner product
$\hproduct_\sim$ is defined as in \eqref{e:sim}.
Moreover, $A_{0,\pm}$ and $V_\pm$ from \eqref{e:AJ}
are selfadjoint in $(\calH_\pm,\hproduct_\sim) = (\calH_\pm,\pm\product)$, $\sigma(A_{0,\pm})\subset\R^\pm\cup\{0\}$ and \eqref{e:sim_adjoint} holds.
We apply Lemma \ref{l:veselic} to the non-negative operators $\pm A_{0,\pm}$ and obtain
\begin{equation*}
\pm I_{\calH_\pm} = \frac 1 \pi\;\,\slim{n\to\infty}\!\int_{1/n}^n\,\left((A_{0,\pm} + it)^{-1} + (A_{0,\pm} - it)^{-1}\right)\,dt.
\end{equation*}
Therefore,
\begin{equation}\label{e:wtJ_ves}
\widetilde J = \frac 1 \pi\;\,\slim{n\to\infty}\!\int_{1/n}^n\,\left((A_{0} + it)^{-1} + (A_{0} - it)^{-1}\right)\,dt,
\end{equation}
where the strong limit is with respect to the norm $\|\cdot\|_\sim := \hproduct_\sim^{1/2}$. But as $\|\cdot\|_\sim$ and $\|\cdot\|$
are equivalent norms, the limit on the right hand side of \eqref{e:wtJ_ves} also exists with respect to the norm $\|\cdot\|$, and
coincides with $\wt J$. The uniform boundedness principle now yields that $\tau_0$ in \eqref{tau0} is finite and that
\begin{equation}\label{e:Jtau}
\|\wt J\|\,\le\,\tau_0<\infty.
\end{equation}
Define the constants
\begin{equation}\label{e:consts}
\nu := \|B\|_\sim,\quad a := \min\sigma(V_+)\quad\text{and}\quad b := \max\sigma(V_-).
\end{equation}
Then $A_{0,+} + V_+$ is bounded from below by $a$ and $\sigma(A_{0,+} + V_+)\subset [a,\infty)$, and
$A_{0,-} + V_-$ is bounded from above by $b$ and $\sigma(A_{0,-} + V_-)\subset (-\infty,b]$. Making use of the
representation \eqref{e:om} of $A = A_0 + V$ and Theorem \ref{l:lmm} we conclude:
\begin{enumerate}
\item[(a)]  $\sigma(A)\backslash\R\subset K_\nu((-\infty,b])\cap K_\nu([a,\infty))$.
\item[(b)]  $(b + \nu,\infty)$ is of positive type with respect to $A$.
\item[(c)] $(-\infty,a - \nu)$ is of negative type with respect to $A$.
\item[(d)]  The growth of the resolvent of $A$ at $\infty$ of order one, and hence of finite order.
\end{enumerate}
The statement in Theorem~\ref{t:main1}~(ii) follows from Theorem~\ref{t:eqthm} if we show that the constants
$\nu$, $a$ and $b$ in \eqref{e:consts} satisfy
$$
\nu\le r, \quad a\ge -d, \quad b\le d,
$$
where $r$ and $d$ are as in \eqref{e:estimates1},
 because in that case the properties (a)-(d) above hold with $\nu$, $a$ and $b$ replaced by $r$, $-d$ and $d$, respectively, and
$$
K_\nu((-\infty,d])\cap K_\nu([-d,\infty))= K_\nu((-d,d)).
$$

First we check $\nu\le r$.
Denote by $\hproduct$ the scalar product corresponding to the Hilbert space norm $\|\cdot\|$. Using the well-known fact that the spectrum of a bounded operator is always a subset of the closure of its numerical range we obtain
\begin{align}
\begin{split}\label{e:rel_proof}
\nu^2
&= \|B\|_\sim^2 = \|B^{\tilde *}B\|_\sim = \sup\{|\la| : \la\in\sigma(B^{\tilde *}B)\}\\
&\le \sup\{|(B^{\tilde *}Bf,f)| : \|f\|=1\}\,\le\,\|B^{\tilde *}\|\,\|B\|\\
&= \|E_-(V|\calH_+)\|\,\|E_+(V|\calH_-)\|\,\le\,\|E_+\|\,\|E_-\|\,\|V\|^2.
\end{split}
\end{align}
As $E_\pm = \frac 1 2(I\pm\wt J)$, it follows from \eqref{e:Jtau} that $\|E_\pm\|\le\frac{1+\tau_0}{2}$ which shows $\nu\le r$.

Next we verify $a\ge -d$.
Denote by $P_\pm$ the orthogonal projection onto $\calH_\pm$ with respect to $\hproduct$. Then, according to Lemma \ref{l:theoretic},
the operator $P_+(J|\calH_+)$ is selfadjoint and uniformly positive in the Hilbert space $(\calH_+,\hproduct)$.
Hence, the same holds for its inverse $E_+(J|\calH_+)$; cf.\ Lemma \ref{l:theoretic}. This implies
\begin{align*}
\min\sigma(P_+(J|\calH_+))
&= \big(\max\sigma(E_+(J|\calH_+))\big)^{-1} = \|E_+(J|\calH_+)\|^{-1}\\
&\ge \|E_+J\|^{-1} \ge \|E_+\|^{-1}\,\ge\,\frac{2}{1+\tau_0}\,.
\end{align*}
Hence, for $f_+\in\calH_+$ we have
$$
[f_+,f_+] = (Jf_+,f_+) = (P_+(J|\calH_+)f_+,f_+)\,\ge\,\frac{2}{1+\tau_0}\,\|f_+\|^2.
$$
We also use
$$(V_+f_+,f_+)_\sim = [V_+f_+,f_+] = [V_+f_+ - B^{\tilde *}f_+,f_+] = [Vf_+,f_+] = (JVf_+,f_+)$$
for $f\in\calH_+$ to conclude
\begin{align*}
a=\min\sigma(V_+)
&=\inf\{(V_+f_+,f_+)_\sim : \|f_+\|_\sim = 1,\,f_+\in\calH_+\}\\
&\ge\inf\{(V_+f_+,f_+)_\sim : \|f_+\|_\sim\le 1,\,f_+\in\calH_+\}\\
&= \inf\{(JVf_+,f_+) : [f_+,f_+]\le 1,\,f_+\in\calH_+\}\\
&\geq\inf\bigl\{(JVf_+,f_+) : \|f_+\|^2\le\tfrac{1+\tau_0}{2},\,f_+\in\calH_+\bigr\}\\
&\geq\inf\bigl\{(JVf,f) : \|f\|^2\le\tfrac{1+\tau_0}{2},\,f\in\calH\bigr\}.
\end{align*}
By assumption $V$ is not non-negative in $(\calH,\product)$ and hence
$$
0>\inf\bigl\{(JVf,f) : \|f\|^2\le\tfrac{1+\tau_0}{2},\,f\in\calH\bigr\}=\inf\bigl\{(JVf,f) : \|f\|^2=\tfrac{1+\tau_0}{2},\,f\in\calH\bigr\}.
$$
Therefore we finally obtain
\begin{align*}
a=\min\sigma(V_+)
&\ge \inf\bigl\{(JVf,f) : \|f\|^2=\tfrac{1+\tau_0}{2},\,f\in\calH\bigr\}\\
&= \frac{1+\tau_0}{2}\,\inf\{(JVf,f) : \|f\| = 1,\,f\in\calH\}\\
&= \frac{1+\tau_0}{2}\,\min\sigma(JV)=-d.
\end{align*}
The proof of the inequality $b\le d$ is analogous and is left to the reader.
This completes the proof of (ii).

As $\infty$ is not a singular critical point of $A_0$, the last
statement in Theorem~\ref{t:main1} on the critical point $\infty$ of $A_0+V$ follows in both cases (i) and (ii)  from
Proposition~\ref{p:curgus} and $\dom A=\dom (A_0+V)$.\qed

\subsection{Proof of Theorem~\ref{t:main2}}\label{ss:proof_main2}
Let $\gamma$ and $\eta$ be as in Theorem~\ref{t:main2}.
Then the operator $A_0 := A + \eta G^{-1}$ is uniformly positive in the Krein space $(\calH,\product)$, $\product = (G\cdot,\cdot)$,
since for $f\in\dom A$ we have
$$
[A_0 f,f] = [Af,f] + \eta[G^{-1}f,f] \ge -\gamma\|f\|^2 + \eta\|f\|^2 = (\eta - \gamma)\|f\|^2.
$$
We will show that $A$ is non-negative over $\ol\C\setminus K_r((-r,r))$.
In particular, $0\in\rho(A_0)$ is not a singular critical point of $A_0$
 and, as $\tau_\eta < \infty$, the point $\infty$ is not a
a singular critical point of $A_0$, see Remark \ref{Mai2012}.
With $V := -\eta G^{-1}$ we have 
$$A = A_0 + V\qquad\text{and}\qquad\Vert V\Vert=\eta\Vert G^{-1}\Vert.$$
As in the proof of Theorem~\ref{t:main1} let $E$ be the spectral function of the operator $A_0$. Then the fundamental
decomposition $\calH = E_+\calH\,[\ds]\,E_-\calH$ in \eqref{e:dec_expl} generates matrix representations
of $A_0$, $V$ and the corresponding fundamental symmetry $\wt J$ as in \eqref{e:AJ} and \eqref{e:V}.
The operator $A_0$ is selfadjoint in the Hilbert space
$(\calH,\hproduct_\sim)$, where $\hproduct_\sim$ is defined in \eqref{e:sim}. With the help of Lemma \ref{l:veselic} we obtain
\begin{align*}
\wt J
&= \frac 1 \pi\;\,\slim{n\to\infty}\!\int_0^n\,\big((A_0 + it)^{-1} + (A_0 - it)^{-1}\big)\,dt\\
&= \frac 1 \pi\;\,\slim{n\to\infty}\!\int_{-n}^n\,(A + \eta G^{-1} - it)^{-1}\,dt\,.
\end{align*}
This yields
$$
\|\wt J\|\,\le\,\tau_\eta.
$$
Following the lines of the proof of Theorem \ref{t:main1}, it suffices to show the inequalities
$$
\nu\le\frac{1+\tau_\eta}{2}\,\|V\|,\quad \min\sigma(V_+)\ge -\frac{1+\tau_\eta}{2}\,\|V\|\quad\text{and}\quad
\max\sigma(V_-)\le \frac{1+\tau_\eta}{2}\,\|V\|,
$$
where $\nu=\Vert E_\mp(V\vert\cal H_\pm) \Vert_\sim$.
The first relation is proved as in \eqref{e:rel_proof} and the remaining two inequalities follow from the selfadjointness of $V_\pm$ in $(\calH_\pm,\hproduct_\sim)$ and
\begin{align*}
\max\{|\min\sigma(V_\pm)|,|\max\sigma(V_\pm)|\}
&= \|V_\pm\|_\sim = \sup\{|\la| : \la\in\sigma(V_\pm)\}\\
&\le \sup\{|(V_\pm f_\pm,f_\pm)| : \|f_\pm\|=1,\,f\in\calH_\pm\}\\
&\le \|V_\pm\| = \|E_\pm(V|\calH_\pm)\|\,\le\,\frac{1+\tau_\eta}{2}\,\|V\|.
\end{align*}
Hence, $A$ is non-negative over $\ol\C\setminus K_r((-r,r))$. As
$\dom A= \dom A_0$, we conclude from Proposition \ref{p:curgus} that $\infty$ is also not a critical point of $A$.
This completes the proof of Theorem~\ref{t:main2}.\qed

\begin{remark}
If, in addition to the assumptions in Theorem \ref{t:main2}, the point $0$ is neither an eigenvalue nor a singular critical point
of the non-negative operator $A_0 := A + \gamma G^{-1}$, then the operator $A$ is non-negative over the set $\ol\C\setminus K_r((-r,r))$, where
\begin{equation}\label{Juni2012}
r = \gamma\,\frac{1+\tau_\gamma}{2}\,\|G^{-1}\|
\end{equation}
and
$$
\tau_\gamma := \frac 1 \pi\,\limsup_{n\to\infty}\,\left\|\int_{1/n}^n\,\big((A_0 + it)^{-1} + (A_0 - it)^{-1}\big)\,dt\right\|.
$$
Indeed, under the assumptions of Theorem \ref{t:main2} for any $\veps > 0$ the operator $A + (\gamma + \veps)G^{-1}$
is selfadjoint in some Hilbert space, and thus the resolvent set of
 $A_0$ is non-empty as $A + (\gamma + \veps)G^{-1}$ and $A_0$ differ just by a bounded operator.
 Therefore, $A_0$ is a non-negative operator in the Krein space
 $(\cal H,\product)$. Moreover,
the point $\infty$ is not a singular critical point of the
operator $A + (\gamma + \veps)G^{-1}$, see Remark \ref{Mai2012},
and by Proposition \ref{p:curgus} $\infty$ is not a singular critical point of $A_0$.
Therefore, the operator $A_0$ admits a matrix representation as in \eqref{e:AJ} and the corresponding fundamental symmetry $\wt J$
satisfies  $
\|\wt J\|\,\le\,\tau_\gamma$ (see the proof of Theorem \ref{t:main2}).
Then by a similar reasoning as in the proof of
Theorem \ref{t:main2} we see that $A$ is non-negative
over  $\ol\C\setminus K_r((-r,r))$, where $r$ is as in \eqref{Juni2012}.
\end{remark}

\section{Differential operators with indefinite weights}\label{s:app}
In this section we apply the results from the previous section to ordinary and partial differential operators with indefinite weights.

\subsection{Indefinite Sturm-Liouville operators}\label{ss:sl}
We consider Sturm-Liouville differential expressions of the form
$$
\calL(f)(x) = \sgn(x)\big(-f''(x) + q(x)f(x)\big),\quad x\in\R\,,
$$
with a real-valued potential $q\in L^\infty(\R)$ and the indefinite weight function $\sgn(\cdot)$.
The corresponding differential operator in $L^2(\R)$ is defined by
$$
Af := \calL(f), \qquad \dom A := H^2(\R).
$$
Here $H^2(\R)$ stands for the usual $L^2$-based second order Sobolov space.
By $\hproduct$ and $\|\cdot\|$ we denote the usual scalar product and its corresponding norm in $L^2(\R)$.
Let $J$ be the operator of multiplication with the function $\sgn(x)$. This operator is obviously selfadjoint
and unitary in $L^2(\R)$. Since the {\it definite} Sturm-Liouville operator
$$
Tf := JAf = -f'' + qf,\quad \dom T := H^2(\R),
$$
is selfadjoint in $L^2(\R)$, it follows that the {\it indefinite} Sturm Liouville operator $A$ is selfadjoint in the Krein space
$(L^2(\R),\product)$, where
$$
[f,g] := (Jf,g) = \int_\R\,f(x)\overline{g(x)}\sgn(x)\,dx,\quad f,g\in L^2(\R).
$$
It is known that the operator $A$ is non-negative over some neighborhood of $\infty$,
but no explicit bounds on the size of this neighborhood exist in the literature.
We first recall a theorem on the qualitative spectral properties of $A$ which can be found
in a slightly different form in  \cite{B07-a,bp,KT09}.
For this, denote by $m_\pm$ the minimum of the essential
spectrum of the selfadjoint operator $T_\pm f_\pm := - f_\pm'' + q_\pm f_\pm$ in $L^2(\R^\pm)$ defined on
$\dom T_\pm := \{f_\pm\in H^2(\R^\pm) : f_\pm(0) = 0\}$, where $q_\pm$
is the restriction of the function $q$ to $\R^\pm$.
The quantities $m_+$ and $m_-$ can be expressed in terms of the potential $q$; if, e.g., $q$ admits limits at $\pm\infty$,
then
\begin{equation*}\label{mpm}
m_\pm=\lim_{x\to\pm\infty}\,q(x).
\end{equation*}

\begin{theorem}
The essential spectrum of the indefinite Sturm-Liouville operator $A$ is the union of the essential spectra of $T_+$ and $-T_-$,
the non-real spectrum of $A$ is bounded and consists of isolated eigenvalues
with finite algebraic multiplicity. Furthermore the following holds.
\begin{enumerate}
\item[{\rm (i)}] If $m_+ > -m_-$ then $\sigma(A)\setminus\R$ is finite;
\item[{\rm (ii)}] If $m_+\leq -m_-$ then $\sigma(A)\setminus\R$ may only  accumulate to points in $[m_+,-m_-]$.
\end{enumerate}
\end{theorem}

To the best of our knowledge explicit bounds on the
non-real spectrum of  indefinite Sturm-Liouville operators
in terms of the potential
$q$ do not exist in the literature. One may expect
that there is a relationship between the maximal magnitude of the non-real eigenvalues of $A$ and the lower bound of the selfadjoint
operator $T = JA$, see, e.g., the numerical examples in \cite{BKT09} and the conjecture in \cite[Remark 11.4.1]{z}.
The following theorem confirms this conjecture and provides
explicit bounds on the non-real eigenvalues of $A$ in terms of the potential $q$. Here only the nontrivial case $\essinf q<0$ is treated
which corresponds to Theorem~\ref{t:main1}~(ii).

\begin{theorem}\label{t:SL2}
Assume that $\essinf q<0$ holds.
Then the indefinite Sturm-Liouville operator $A$ is non-negative over $\ol\C\setminus K_r((-d,d))$, where
$$
r := 5\|q\|_\infty,\quad\text{and}\quad d := -5\,\essinf_{x\in\R}q(x)>0.
$$
The non-real spectrum of $A$ is contained in
$$
 K_r((-d,d))\cap\{\la\in\C : |\Im\la|\le 2\|q\|_\infty\}.
$$
\end{theorem}
\begin{proof}
The proof of Theorem~\ref{t:SL2} is split into 4 parts. The first step is preparatory and connects the present problem
with the abstract setting in Theorem~\ref{t:main1}. In the second part a Krein type resolvent formula is provided
which is essential for the main estimates in the last two parts of the proof.

\
\\
\noindent{\it 1.\ Preparation: } The operator $A_0$,
\begin{equation*}
A_0 f:=-\sgn(\cdot)f^{\prime\prime},\qquad \dom A_0=H^2(\R),
\end{equation*}
is selfadjoint and non-negative in the Krein space $(L^2(\R),[\cdot,\cdot])$. Furthermore, we have $\sigma(A_0) = \R$,
and neither $0$ nor $\infty$ is a singular critical point of $A_0$, see \cite{cn}. Hence $A_0$ satisfies the
assumptions in Theorem~\ref{t:main1}. Define $V$ as
\begin{equation*}
 Vf:=\sgn(\cdot) qf,\qquad \dom V=L^2(\R),
\end{equation*}
so that $V$ is a bounded selfadjoint operator in $(L^2(\R),[\cdot,\cdot])$ with
$$
\min\sigma(JV) = \essinf_{x\in\R}\,q(x)\quad\text{and}\quad\Vert V\Vert = \Vert q\Vert_\infty.
$$
By Theorem \ref{t:main1} the operator
$$
A = A_0 + V
$$
is non-negative over $\ol\C\setminus  K_r((-d,d))$, where
$$
r = \frac{1+\tau_0}{2}\,\|q\|_\infty,\qquad d = -\frac{1+\tau_0}{2}\,\essinf_{x\in\R}q(x),
$$
and
$$
\tau_0 = \frac 1 \pi\,\limsup_{n\to\infty}\,\left\|\int_{1/n}^n\,\left((A_0 + it)^{-1} + (A_0 - it)^{-1}\right)\,dt\,\right\|.
$$
Therefore, we have to show
\begin{equation}\label{e:assert1}
\tau_0\le 9
\end{equation}
and
\begin{equation}\label{e:assert2}
\sigma(A)\setminus\R\,\subset\,\{\la\in\C : |\Im\la|\le 2\|q\|_\infty\}.
\end{equation}

\ \\
\noindent{\it 2.\ Krein's resolvent formula: }In the following the resolvent of $A_0$ will be expressed via a Krein type resolvent
formula in terms of the resolvent of the diagonal operator matrix
\begin{equation*}
B_0:=\begin{pmatrix} B_+ & 0 \\ 0 & B_- \end{pmatrix}\quad \text{in}\quad L^2(\R)=L^2(\R^+)\times L^2(\R^-),
\end{equation*}
where $B_\pm f_\pm:=\mp f_\pm^{\prime\prime}$ are the selfadjoint Dirichlet operators in $L^2(\R^\pm)$ which are defined
on $\dom B_\pm = \{f_\pm \in H^2(\R^\pm):f_\pm(0)=0\}$. Here and in the following we will denote the restrictions of a function $f$
defined on $\R$ to $\R^\pm$ by $f_\pm$. For $\la = re^{it}$, $r > 0$, $t\in [0,2\pi)$, we set $\sqrt{\la} := \sqrt{r}e^{it/2}$.
For $\la\in\C\setminus\R$ define the function
$$
f_\lambda(x) :=
\begin{cases}
e^{i\sqrt{\lambda}x}, &x> 0,\\
e^{-i\sqrt{-\lambda}x}, &x < 0.
\end{cases}
$$
The function $f_\lambda$ is a solution of the equations $\mp f_\pm^{\prime\prime} = \lambda f_\pm$ in $L^2(\R^\pm)$. From $\overline{\sqrt{\lambda}} = -\sqrt{\overline\lambda}$ it follows that $f_{\ol\la} = \ol{f_\la}$ holds for all $\la\in\C\setminus\R$. Let us now prove that for all $f\in L^2(\R)$ and $\lambda\in \C\setminus\R$ we have
\begin{equation}\label{e:krf}
(A_0 - \la)^{-1}f = (B_0 - \la)^{-1}f - \frac{[f,f_{\ol\la}]}{i(\sqrt{\la}+\sqrt{-\la})}f_\la.
\end{equation}
Since $\sigma(A_0)=\R$ and $B_0$ is selfadjoint in the Hilbert space $L^2(\R)$
it follows that the resolvents of $A_0$ and $B_0$ in \eqref{e:krf} are defined for all $\lambda\in \C\setminus\R$.
In particular, for $\lambda\in \C\setminus\R$ there exists $g\in\dom B_0$ such that $f=(B_0-\lambda)g$ holds.
The right hand side in \eqref{e:krf} has the form
\begin{equation}\label{hhh}
h:=g-\frac{[(B_0-\lambda)g,f_{\ol\lambda}]}{i(\sqrt{\la}+\sqrt{-\la})}f_\la
\end{equation}
and we will show that the function $h$ belongs to $\dom A_0=H^2(\R)$. As $g_\pm$ and $f_{\lambda,\pm}$ are elements of $H^2(\R^\pm)$
the same is true for $h_\pm$. Moreover, $g$ and $f_\lambda$ are continuous at $0$ and so is $h$.
Hence it remains to
check that $h^\prime$ is continuous at $0$.
For this note first that
\begin{equation*}\label{qwert}
[(B_0-\lambda)g,f_{\ol\lambda}]=\int_0^\infty (-g_+^{\prime\prime}-\lambda g_+)\overline{f_{\ol\lambda,+}}\, dx
- \int_{-\infty}^0 (g_-^{\prime\prime}-\lambda g_-)\overline{f_{\ol\lambda,-}}\, dx
\end{equation*}
and since $f_{\ol\lambda,\pm}$ are solutions of $\mp f_\pm^{\prime\prime}=\ol\lambda f_\pm$, integration by parts yields
\begin{equation*}[(B_0-\lambda)g,f_{\ol\lambda}]=
g_+^\prime(0)\overline{f_{\ol\lambda,+}(0)}-g_+(0) \overline{f^\prime_{\ol\lambda,+}(0)}-
g_-^\prime(0)\overline{f_{\ol\lambda,-}(0)}+g_-(0) \overline{f^\prime_{\ol\lambda,-}(0)}.
\end{equation*}
As $g\in\dom B_0$ we have $g_\pm(0)=0$ and together with $f_{\ol\lambda,\pm}(0)=1$ we find
\begin{equation*}
 [(B_0-\lambda)g,f_{\ol\lambda}]= g_+^\prime(0)-g^\prime_-(0).
\end{equation*}
Therefore we obtain for the derivatives $h_\pm^\prime$ on $\R^\pm$ of the function $h$ from \eqref{hhh}:
\begin{equation*}
h^\prime_\pm=g^\prime_\pm-\frac{g_+^\prime(0)-g^\prime_-(0)}{i(\sqrt{\la}+\sqrt{-\la})}f_{\lambda,\pm}^\prime
\end{equation*}
and as $f^\prime_{\lambda,+}(0)=i\sqrt{\lambda}$ and $f^\prime_{\lambda,-}(0)=-i\sqrt{-\lambda}$ we conclude
\begin{equation*}
h^\prime_+(0)-h^\prime_-(0)= \bigl(g_+^\prime(0)-g^\prime_-(0)\bigr)-
\frac{g_+^\prime(0)-g^\prime_-(0)}{i(\sqrt{\la}+\sqrt{-\la})}\bigl(f^\prime_{\lambda,+}(0)-f^\prime_{\lambda,-}(0)\bigr)=0,
\end{equation*}
that is, $h^\prime$ is continuous at $0$ and therefore $h\in\dom A_0$.

Now a straightforward computation shows that $(A_0-\lambda)h=(B_0-\lambda)g=f$ holds and hence
the resolvent of $A_0$ is given by \eqref{e:krf}.

\ \\
\noindent{\it 3.\ Proof of {\rm\eqref{e:assert1}}: }Let $f\in L^2(\R)$ and $t > 0$. Then \eqref{e:krf} and $\sqrt{it} + \sqrt{-it} = i\sqrt{2t}$ yield
\begin{align}
\begin{split}\label{e:needed_later}
[(A_0 + it)^{-1}f + (A_0 - it)^{-1}f,f] = \,
& [(B_0 + it)^{-1}f + (B_0 - it)^{-1}f,f]\\
& + \frac 2 {\sqrt{2t}}\,\Re([f,f_{it}][f_{-it},f]).
\end{split}
\end{align}
With $g(x) := |f(x)| + |f(-x)|$, $x\in\R^+$, we have
$$
\big|[f,f_{it}][f_{-it},f]\big|\,\le\,\left(\int_0^\infty\,e^{-x\sqrt{t/2}}g(x)\,dx\right)^2
$$
and thus for $n\in\N$
\begin{align*}
\int_{1/n}^n\,\frac{|[f,f_{it}][f_{-it},f]|}{\sqrt{2t}}\,dt
&\le \int_{1/n}^n\,\int_0^\infty\,\int_0^\infty\,g(x)g(y)\,\frac{e^{-(x+y)\sqrt{t/2}}}{\sqrt{2t}}\,dy\,dx\,dt\\
&= 2\,\int_0^\infty\,\int_0^\infty\,\frac{g(x)g(y)}{x+y}\left(e^{-\frac{x+y}{\sqrt{2n}}} - e^{-\frac{\sqrt{n}(x+y)}{\sqrt{2}}}\right)\,dy\,dx\\
&\le 2\,\int_0^\infty\,\int_0^\infty\,\frac{g(x)g(y)}{x+y}\,dy\,dx.
\end{align*}
From Hilbert's inequality (see, e.g., \cite{hlp}) it follows that
\begin{equation*}
\int_{1/n}^n\,\frac{|[f,f_{it}][f_{-it},f]|}{\sqrt{2t}}\,dt\,\le\,2\pi\|g\|_{L^2(\R^+)}^2\,\le\,4\pi\|f\|^2.
\end{equation*}
and therefore we have
\begin{equation}\label{e:2nd}
\frac{1}{\pi}\int_{1/n}^n\,\frac 2 {\sqrt{2t}}\,\Re([f,f_{it}][f_{-it},f]) \,dt\,\le\,8\|f\|^2.
\end{equation}

Denote by $E_0$ ($E_\pm$) the spectral function of the operator $A_0$ ($B_\pm$, respectively). Moreover define
$$
J_{n,\pm} := \frac 1 \pi\,\int_{1/n}^n\,\left((B_\pm + it)^{-1} + (B_\pm - it)^{-1}\right)\,dt,
$$
and, in addition,
$$
J_{n,0} := \frac 1 \pi\,\int_{1/n}^n\,\left((A_0 + it)^{-1} + (A_0 - it)^{-1}\right)\,dt.
$$
As $\pm B_\pm$ is non-negative, it follows that $J_{n,\pm}$ converges strongly to $\pm I_{L^2(\R^\pm)}$ when $n\rightarrow\infty$, 
cf.\ Lemma \ref{l:veselic}. Moreover, by
$$(B_\pm + it)^{-1} + (B_\pm - it)^{-1}= 2B_\pm (B_\pm + it)^{-1} (B_\pm - it)^{-1}$$
the sequence $\pm J_{n,\pm}$ is an increasing sequence of non-negative selfadjoint operators in $L^2(\R^\pm)$. This implies
\begin{align*}
\frac 1 \pi\,\int_{1/n}^n\,[(B_0 + it)^{-1}f + (B_0 - it)^{-1}f,f]\,dt
&= (J_{n,+}f_+,f_+) - (J_{n,-}f_-,f_-)\\
&\le \|f_+\|_{L^2(\R^+)}^2 + \|f_-\|_{L^2(\R^-)}^2 = \|f\|^2.
\end{align*}
Now, from \eqref{e:needed_later} and \eqref{e:2nd} it follows that
$$
[J_{n,0}f,f]\,\le\,\|f\|^2 + 8\|f\|^2 = 9\|f\|^2,\quad f\in L^2(\R).
$$
 Hence, as $J = \sgn(\cdot)$ is unitary in $L^2(\R)$, we finally obtain
$\|J_{n,0}\| = \|J J_{n,0}\|\le 9$
and thus $\tau_0 = \limsup_{n\to\infty}\,\|J_{n,0}\|\le 9$. That is, \eqref{e:assert1} holds.

\ \\
\noindent{\it 4.\ Proof of {\rm\eqref{e:assert2}}: }Let $\la\in\C^+$. Using \eqref{e:krf}, for $f\in L^2(\R)$ we obtain
$$
\|(A_0 - \la)^{-1}f\|\,\le\,\|(B_0 - \la)^{-1}f\| + \frac{|[f,f_{\ol\la}]|\,\|f_\la\|}{|\sqrt{\la} + \sqrt{-\la}|}\,\le\,\frac{\|f\|}{|\Im\la|} + \frac{\|f_\la\|^2\|f\|}{\sqrt{2|\la|}}\,.
$$
Here we have used the identity $f_{\ol\la} = \ol{f_\la}$. Let $r > 0$ and $t\in (0,\pi)$ such that $\la = re^{it}$. Then
\begin{align*}
\|f_\la\|^2
&= \int_0^\infty\,\left|e^{i\sqrt{r}e^{it/2}x}\right|^2\,dx + \int_{-\infty}^0\,\left|e^{-i\sqrt{r}e^{i(t+\pi)/2}x}\right|^2\,dx\\
&= \int_0^\infty\,e^{-2\sqrt{r}\sin(t/2)x}\,dx + \int_{-\infty}^0\,e^{2\sqrt{r}\cos(t/2)x}\,dx\\
&= \frac{\cos(t/2) + \sin(t/2)}{\sqrt{r}\,\sin(t)}
\,\le\,\frac{\sqrt{2}}{\sqrt{r}\,\Im(\la/r)} = \frac{\sqrt{2|\la|}}{|\Im\la|}.
\end{align*}
Therefore,
$$
\|(A_0 - \la)^{-1}\|\,\le\,\frac{2}{|\Im\la|}\,.
$$
The same estimate holds for $\la\in\C^-$. Now, assume that $|\Im\la| > 2\|q\|_\infty=2\| V\|$. Then $\|V(A_0 - \la)^{-1}\| < 1$, and hence
$$
A - \la = A_0 - \la + V = \left(I + V(A_0 - \la)^{-1}\right)(A_0 - \la)
$$
is boundedly invertible, which proves \eqref{e:assert2}.\qed
\end{proof}

\subsection{Second order elliptic operators}\label{ss:pde}
Let $\Omega\subset\R^n$ be a domain and let $\ell$ be the "formally
symmetric" uniformly
elliptic second order differential expression
\begin{equation*}
\ell (f)(x):=-\sum_{j,k=1}^n \left( \frac{\partial}{\partial x_j} a_{jk}
\frac{\partial f}{\partial x_k}\right)(x)+a(x)f(x),\quad x\in\Omega,
\end{equation*}
with bounded coefficients $a_{jk}\in C^\infty(\Omega)$ satisfying
$a_{jk}(x)=\overline{a_{kj}(x)}$
for all $x\in\Omega$ and $j,k=1,\dots,n$, the function $a\in
L^\infty(\Omega)$ is real valued and
\begin{equation*}
\sum_{j,k=1}^n a_{jk}(x)\xi_j\xi_k\geq C\sum_{k=1}^n\xi_k^2
\end{equation*}
holds for some $C>0$, all $\xi=(\xi_1,\dots,\xi_n)^\top\in\R^n$ and
$x\in\Omega$. With the differential expression $\ell$ we associate the
elliptic differential operator
\begin{equation*}
Tf:=\ell (f),\qquad \dom T=\bigl\{f\in H^1_0(\Omega):\ell (f) \in
L^2(\Omega)\bigr\},
\end{equation*}
where $H^1_0(\Omega)$ stands for the closure of $C_0^\infty(\Omega)$ in
the Sobolev space $H^1(\Omega)$.
It is well known that $T$ is an unbounded selfadjoint operator in the
Hilbert space $(L^2(\Omega),(\cdot,\cdot))$ with spectrum
semibounded from below by $\essinf\,a$; cf. \cite{EE}.

Let $w$
be a real valued function such that $w,w^{-1}\in L^\infty(\Omega)$
and each of the sets
\begin{equation*}
 \Omega_+:=\bigl\{x\in\Omega: w(x)>0\bigr\}\quad\text{and}\quad
\Omega_-:=\bigl\{x\in\Omega: w(x)<0\bigr\}
\end{equation*}
has positive Lebesgue measure.
We define a second order elliptic differential
expression $\calL$ with the indefinite weight $w$ by
\begin{equation*}
\calL(f)(x):=\frac{1}{w(x)}\,\ell (f)(x),\qquad x\in\Omega.
\end{equation*}
The multiplication operator $G_w f=w f$, $f\in L^2(\Omega)$, is an
isomorphism in $L^2(\Omega)$
with inverse $G_w^{-1}f=w^{-1}f$, $f\in L^2(\Omega)$, and gives rise to
the Krein space inner product
\begin{equation*}
[f,g]:=(G_w f,g)=\int_\Omega f(x)\overline{g(x)}\,w(x)\,dx,\qquad f,g\in
L^2(\Omega).
\end{equation*}
The differential operator associated with $\calL$ is defined as
\begin{equation}\label{indefa}
Af=\calL (f),\qquad \dom A=\bigl\{f\in H^1_0(\Omega):\calL (f) \in
L^2(\Omega)\bigr\}.
\end{equation}
Since for $f\in H^1_0(\Omega)$ we have $\ell(f)\in L^2(\Omega)$ if and only if $\calL(f)\in L^2(\Omega)$
it follows that  $\dom A=\dom T$
and $A= G_w^{-1}T$ hold. Hence
$A$ is a selfadjoint operator in the Krein space $(L^2(\Omega),[\cdot,\cdot])$.

In order to illustrate Theorem~\ref{t:main2}
for the indefinite elliptic operator $A$ we assume from now on that
$$\min\sess(T)\leq 0$$
 holds. This also implies that the domain $\Omega$ is unbounded
as otherwise $\sess (T)=\emptyset$. A discussion of the cases $\sess
(T)=\emptyset$ and $\min\sess(T)> 0$ is contained in \cite{B11}, see also
\cite{F90,P89}.
Fix some $\eta>0$ such that $-\eta<\min\sigma(T)$ and define the
spaces $\calH_s$, $s\in[0,2]$, as the domains of the
$\frac{s}{2}$-th powers of the uniformly positive operator $T+\eta$ in
$L^2(\Omega)$,
\begin{equation*}
\calH_s:=\dom \bigl((T+\eta)^\frac{s}{2}\bigr),\qquad s\in [0,2].
\end{equation*}
Note that $\calH=\calH_0$, $\dom T=\calH_2$ and the form domain of $T$
is $\calH_1$. The spaces $\calH_s$
become Hilbert spaces when they are equipped with the usual inner
products, the induced
topologies do not depend on the particular choice of $\eta$; cf.\ \cite{k}.

The following theorem is a direct consequence of Theorem~\ref{t:main2} and the
considerations in \cite[Theorem 2.1~(iii)]{CN94} and \cite{cu}; cf. \cite[Theorem 4.4]{B11}.

\begin{theorem}\label{bigthm}
Let $A$ be the indefinite elliptic operator in \eqref{indefa}, $\eta>0$ as above, and
assume that there exists a bounded uniformly positive operator $W$ in
$(L^2(\Omega),[\cdot,\cdot])$
such that $W\calH_s\subset\calH_s$ holds for some $s\in (0,2]$.
Then $A$ is non-negative over $ \ol{\mathbb C} \setminus K_r((-r,r))$, where
$$r=\eta\,\frac{1+\tau_\eta}{2}\,\Vert w^{-1}\Vert_\infty\quad\text{and}\quad
\tau_\eta := \frac 1 \pi\,\limsup_{n\to\infty}\,\left\|\int_{-n}^n\,(A+\eta G_w^{-1} - it)^{-1} \,dt\,\right\|\,<\,\infty.
$$
Moreover, $\infty$ is not a singular critical point of $A$.
\end{theorem}

In the next corollary the special case $\Omega=\R^n$ is treated; cf. \cite[Theorem 5.4]{B11}. From now on we assume that
$\Omega=\R^n$ and $\Omega_\pm=\{x\in\R^n:\pm w(x)>0\}$
consist of finitely many connected components with compact smooth boundaries, and that the coefficients
$a_{jk}\in C^\infty(\R^n)$ and their derivatives
are uniformly continuous and bounded. Note that either $\Omega_+$ or $\Omega_-$ is bounded.

\begin{corollary}\label{cor1}
Suppose that for some $s\in(0,\frac{1}{2})$ the Sobolev spaces $H^s(\Omega_\pm)$ are invariant
under multiplication with $w\!\upharpoonright_{\Omega_\pm}$. Then the indefinite elliptic operator $A$ is non-negative over
$\ol{\mathbb C} \setminus K_r((-r,r))$, where $r$ is as in in Theorem~\ref{bigthm}. Moreover, $\infty$ is not a singular critical point of $A$.
\end{corollary}

A sufficient criterion for the invariance of  $H^s(\Omega_\pm)$ in the above corollary can be deduced from \cite[Theorem 1.4.1.1]{G85}.

\begin{corollary}
Suppose that the weight function $w\!\upharpoonright_{\Omega_\pm}$ belongs to some H\"{o}lder space $C^{0,\alpha}(\Omega_\pm)$, 
$\alpha\in(0,\frac{1}{2})$, and that outside of some bounded set $w$ is equal to a constant. Then the claim in Corollary~\ref{cor1} holds.
\end{corollary}

\begin{proof}{\hspace*{-.1cm}{\it of Corollary~\ref{cor1}}\,}
The assumptions on the coefficients $a_{jk}$ imply
\begin{equation*}
\dom A=\dom T=H^2(\R^n)\qquad\text{and}\qquad \calH_s=H^s(\R^n),\qquad s\in [0,2],
\end{equation*}
by  elliptic regularity and interpolation. Since the Sobolev spaces $H^s(\Omega_\pm)$ are assumed to be invariant
under multiplication with the functions $w\!\upharpoonright_{\Omega_\pm}$ for some $s\in(0,\frac{1}{2})$ it follows from
\cite[Lemma~5.1]{B11} that also $H^s(\R^n)$ is invariant under multiplication with $w$, that is, $G_w\calH_s\subset\calH_s$
holds. Furthermore, $G_w$ is uniformly positive in $(L^2(\Omega),[\cdot,\cdot])$ since
$$[G_w f,f]=(G_w^2 f,f)\geq \essinf w^2 \Vert f\Vert^2$$
and $w^{-1}\in L^\infty(\R^n)$.
Now the assertion follows from Theorem~\ref{bigthm}.\qed
\end{proof}

\section*{Acknowledgements}
The support from the Deutsche For\-schungs\-ge\-mein\-schaft (DFG) under the grants BE 3765/5-1 and TR 903/4-1
is gratefully acknowledged.

\section*{Contact information}
Jussi Behrndt: Institut f\"ur Numerische Mathematik, Technische Universit\"at Graz, Steyrergasse 30, 8010 Graz, Austria, behrndt@tugraz.at

\vspace{0.4cm}\noindent
Friedrich Philipp: Institut f\"ur Mathematik, MA 8-1, Technische Universit\"at Berlin, Stra\ss e des 17.\ Juni 136, 10623 Berlin, Germany, fmphilipp@gmail.com

\vspace{0.4cm}\noindent
Carsten Trunk: Institut f\"ur Mathematik, Technische Universit\"at Ilmenau, Postfach 10 05 65, 98684 Ilmenau, Germany, carsten.trunk@tu-ilmenau.de
\end{document}